\magnification=\magstep1
An automorphic form related to cubic surfaces. First draft. 1997.

R. E. Borcherds. 

(Remark added 2000: See the paper ``Cubic Surfaces and Borcherds Products''
by Allcock and Freitag, math.AG/0002066)

In the paper [A-C-T] the authors showed that the moduli space of cubic
surfaces was $(CH^4\backslash H)/G$, where $H$ is the union of the
reflection hyperplanes of the reflection group $G$. The purpose of
this note is to construct an automorphic form, called the discriminant,
on complex
hyperbolic space $CH^4$
whose zeros are exactly the reflection hyperplanes of $G$,
each with multiplicity 1.

The complex hyperbolic space of $E^{1,4}$ embeds naturally in the
Grassmannian of 2 dimensional positive definite subspaces of
the underlying even integral lattice $M$ of $E^{1,4}$. We will construct
the automorphic form on complex hyperbolic space by constructing
an automorphic form $\Psi$ on the Grassmannian and restricting it to 
complex hyperbolic space. 
The construction of $\Psi$
is similar to that of
example 13.7 of [B]. 
The lattice $M$ is isomorphic to $A_2\oplus A_2^4(-1)$, so it
has dimension 10 and determinant
$3^5$.  We let the elements $e_\gamma$ for $\gamma\in M'/M$
stand for the obvious basis of the group ring
$C[M'/M]$. We will construct a modular form of
weight $-3$ and type $\rho_M$ for $SL_2(Z)$ which is holomorphic
on the upper half plane and meromorphic at cusps,
by which we mean a holomorphic vector valued
function 
$f=\sum_{\gamma\in M'/M}e_\gamma f_\gamma(\tau)$
such that 
$$
f_\gamma(\tau+1)= e^{\pi i \gamma^2}f_\gamma(\tau),
\qquad
f_\gamma(-1/\tau) = -i3^{-5/2}\tau^{-3}
\sum_{\delta\in M'/M} e^{-2\pi i (\delta,\gamma)}
f_\delta(\tau)
$$

Recall that if $\gamma,\delta\in M'/M$ then $\gamma^2$ is well defined mod 2
and $(\gamma,\delta)$ is well defined mod 1. There are 4 orbits
of vectors $\gamma\in M'/M$ under ${\rm Aut}(M)$, which we name as follows:
a nonzero vector $\gamma$ with $\gamma^2/2\equiv n/3 \bmod 1$ will be
called a vector of type $n$ ($n=0,1,2$), and the zero vector will be called
a vector of type $00$. 

We will construct a modular form $f$ such that the components
$f_\gamma$ of $f$ are given by functions $f_{00}$, $f_0$, $f_1$, or $f_2$
depending only on the type of $\gamma$. We now work out the conditions
that these functions have to satisfy for $f$ to transform correctly.
In order to check that $f$ is a modular form under $\tau\mapsto -1/\tau$
we need to know, given some fixed vector $u$, how many vectors $v$
there are with given type and given inner product with $u$. These
numbers are given in the following table. 
$$
\matrix{
\hbox{type of $u$}&00&00&00   &00&0&0&0&0&1&1&1&1&2&2&2&2 \cr
\hbox{type of $v$}&00&0&1&2   &00&0&1&2&00&0&1&2&00&0&1&2 \cr
(u,v)\equiv 0&
1&80&90&72& 1&26&36&18&  1&32&24&24& 1&20&30&30\cr
(u,v)\equiv1/3&
0&0&0&0&    0&27&27&27&  0&24&33&24& 0&30&30&21\cr 
(u,v)\equiv2/3&
0&0&0&0&    0&27&27&27&  0&24&33&24& 0&30&30&21\cr 
} $$

Using this table we see that  $f_{00}$, $f_0$, $f_1$, 
and $f_2$ have to satisfy the equations
$$
\eqalign{
&f_{00}(\tau+1)=f_{00}(\tau), \quad
f_{00}(-1/\tau)
=-i3^{-5/2}\tau^{-3}(f_{00}(\tau)+80f_0(\tau)+90f_1(\tau)+72f_2(\tau))\cr
&f_{0}(\tau+1)=f_{0}(\tau),\quad
f_{0}(-1/\tau)
=-i3^{-5/2}\tau^{-3}(f_{00}(\tau)-f_0(\tau)+9f_1(\tau)-9f_2(\tau))\cr
&f_{1}(\tau+1)=e^{2\pi i/3}f_{1}(\tau),\quad
f_{1}(-1/\tau)
=-i3^{-5/2}\tau^{-3}(f_{00}(\tau)+8f_0(\tau)-9f_1(\tau))\cr
&f_{2}(\tau+1)=e^{4\pi i/3}f_{2}(\tau),\quad
f_{2}(-1/\tau)
=-i3^{-5/2}\tau^{-3}(f_{00}(\tau)-10f_0(\tau)+9f_2(\tau))\cr
}$$

One solution of these equations is given as follows.  
$$
\eqalign{
f_{00}(\tau)&=
24\eta(3\tau)^3\eta(\tau)^{-9}=24(1+9q+54q^2+O(q^3)) \cr
f_0(\tau)&=-3\eta(3\tau)^3\eta(\tau)^{-9}=-3+O(q) \cr
f_1(\tau) &=0\cr
f_2(\tau)&=
\eta(\tau/3)^3\eta(\tau)^{-9}+3\eta(3\tau)^3\eta(\tau)^{-9}
=q^{-1/3}+14q^{2/3}+92q^{5/3}+O(q^{8/3})\cr
}
$$
Most of the transformations follow formally from the functional
equations $\eta(\tau+1)=e^{2\pi i /24}\eta(\tau)$ and
$\eta(-1/\tau)=\sqrt{\tau/i}\eta(\tau)$ of $\eta$.  The only one which
takes slightly more work is the transformation of $f_2$ under
$\tau\mapsto\tau+1$, and this follows from the identity
$\eta(\tau)^3=\sum_{n\in Z}(4n+1)q^{(4n+1)^2/8}$ and its consequence
$$\eta(\tau/3)^3\eta(\tau)^{-9}+\eta((\tau+1)/3)^3\eta(\tau+1)^{-9}
+\eta((\tau+2)/3)^3\eta(\tau+2)^{-9} = -9\eta(3\tau)^3\eta(\tau)^{-9}.$$

By theorem 13.3 of [B] there is an automorphic form $\Psi$
on the symmetric space
of $M$ with the
following properties. It has weight 12 = (coefficient of $q^0$ in
$f_{00}$)/2.  The zeros of $\Psi$
correspond to the negative powers of $q$ in $f$,
so are zeros of order 1 orthogonal to all the norm $-2/3$ vectors of
$M'$.  $\Psi$ is holomorphic on the symmetric space, and therefore
holomorphic at cusps as well by the Koecher boundedness principle. $\Psi$
is an automorphic form for some one dimensional representation of 
${\rm Aut}(M)$.

By restricting $\Psi$ to complex hyperbolic space we get an
automorphic form which has zeros of order 3 along the reflection
hyperplanes (because $f$ has zeros of order 1, but every reflection
hyperplane is the restriction of 3 hyperplanes of the symmetric space
of $M$).  So by taking the cube root of the restriction of $\Psi$ we get an
automorphic form for a one dimensional character of $G$ whose zeros
are exactly the reflection hyperplanes with multiplicity 1.

References.

[A-C-T] D. J. Allcock, J. Carlson, D. Toledo,
A Complex Hyperbolic Structure for Moduli of Cubic Surfaces, 
alg-geom/9709016

[B] R. E. Borcherds, Automorphic forms with singularities on Grassmannians,
alg-geom/9609022, to appear in Invent. Math. 
\bye